# SINGULAR NORMAL FORM
# FOR THE PAINLEVÉ EQUATION P1

OVIDIU COSTIN AND RODICA D. COSTIN

ABSTRACT. We show that there exists a rational change of coordinates of Painlevé's P1 equation $y'' = 6y^2 + x$ and of the elliptic equation $y'' = 6y^2$ after which these two equations become analytically equivalent in a region in the complex phase space where $y$ and $y'$ are unbounded. The region of equivalence comprises all singularities of solutions of P1 (i.e. outside the region of equivalence, solutions are analytic). The Painlevé property of P1 (that the only movable singularities are poles) follows as a corollary. Conversely, we argue that the Painlevé property is crucial in reducing P1, in a singular regime, to an equation integrable by quadratures.

## 1. INTRODUCTION

The problem of determining which classes of nonlinear differential equations can define new transcendents (special functions having good properties), received a special attention in the last century, especially due to the emphasis on finding "explicit" solutions to differential equations. Fuchs had the intuition that the appropriate condition these equations must satisfy is that their solutions have no movable branch points. This feature of an equation is now known as the Painlevé property and proved to be a very relevant characteristic, in a wide range of problems. Fuchs' study was pursued by Briot and Bouquet, and then by Painlevé [10] and Gambier who showed that there are no new transcendents coming from first order equations, but there are six second order equations which define new special functions. These equations (now denoted usually as P1 to P6) were discovered as a result of a purely theoretical quest, but they later arose naturally in many distinct physical applications (see, e.g., [4] and [7]). Linearization of second order Painlevé equations through the isomonodromic transformation method [2], [4], [8], [9] is one of the most important recent developments. To this date, higher order equations have not yet been classified from the point of view of Painlevé integrability.

Perhaps surprisingly, proving the Painlevé property of an equation turns out to be quite difficult (although if one *assumes* that singularities are described locally by convergent power-logarithmic series, then it is usually easy to check for the absence of movable branch points) and some of the classical proofs for Painlevé equations have been subsequently challenged. See also [5] and [6].

The Painlevé property, being shared by all solutions, must reveal a particular structure of the equation itself. We show in fact that P1 is equivalent with an equation integrable by quadratures, in a region in the phase space where the solutions $y(x)$

Research at MSRI was supported in part by NSF grants DMS-9701755 and DMS-9704968.

are singular. The equivalence also has the implication that the solutions of P1 are meromorphic, and is a natural and rigorous way to prove the Painlevé property.

We expect our technique to work for other equations having the Painlevé property, as well.

The existence of a simple integrable singular normal form of P1 is tied to the special integrability properties of this equation, and is not merely a consequence of the fact that (2.1) "approximately equals" $y'' = 6y^2$ when $y, y', y''$ are large. Approximations near singularities are usually very *unstable*. For instance, the modification $Y'' = 6Y^2 + x^2$ of P1 also seems to be approximated by the equation $y'' = 6y^2$ for large $Y, Y''$, especially if $x$ is small, but these two equations are *not* equivalent in mentioned regime. The obstruction in the equivalence is the presence of "bad" logarithmic terms in the Frobenius series of $Y(x)$ near its movable singularities, and the fact that their *absence* in $y(x)$ is stable under combinations of analytic and rational transformations.

## 2. Main Results

We consider the Painlevé P1 equation

$$(2.1) \qquad \frac{d^2 y}{dx^2} = 6y^2 + x$$

and show that there exists a transformation of the independent variable only (i.e., of the form $x = F(t, u, u')$, $y = u$) which is an equivalence of (2.1) to the elliptic equation

$$(2.2) \qquad \frac{d^2 u}{dt^2} = 6u^2$$

in regions of the phase space where the dependent variables are large.

The regularity of the transformation giving the equivalence appears more clearly after making a rational transformation of the dependent variables of (2.1) and (2.2)

$$(2.3) \qquad v_1 = -y^{-2} \frac{dy}{dx}, \quad w_1 = y^3 \left(\frac{dy}{dx}\right)^{-2}$$

and, respectively

$$(2.4) \qquad v = -u^{-2} \frac{du}{dt}, \quad w = u^3 \left(\frac{du}{dt}\right)^{-2}$$

Then (2.1) and (2.2), written in the new variables,

$$(2.5) \qquad \begin{cases} \frac{dv_1}{dx} = \frac{2}{w_1} - 6 - x w_1^2 v_1^4 \\ \frac{dw_1}{dx} = \frac{12}{v_1} \left(w_1 - \frac{1}{4}\right) + 2x w_1^3 v_1^3 \end{cases}$$

and, respectively

$$
(2.6) \quad \begin{cases} \frac{dv}{dt} = \frac{2}{w} - 6 \\ \\ \frac{dw}{dt} = \frac{12}{v}\left(w - \frac{1}{4}\right) \end{cases}
$$

are *analytically equivalent* on polydisks centered at $(x_0, 0, 1/4)$ ($x_0 \in \mathbb{C}$ arbitrary). We note that the centers of the polydisks correspond to the point at infinity for $y$ and $u$. In remark 5 below we explain how this equivalence translates into properties of the solutions of P1.

**Proposition 1.** *There exists a map $(t, v, w) \mapsto \Gamma(t, v, w) = (x, v_1, w_1)$, in the extended phase space, of the form*

$$
(2.7) \quad \begin{cases} x = \gamma(t, v, w) \\ v_1 = \eta(t, v, w) \\ w_1 = \theta(t, v, w) = v^2 w \eta(t, v, w)^{-2} \end{cases}
$$

*which transforms (2.5) into (2.6) and has the following properties:*

*(i) $\Gamma$ is holomorphic in a neighborhood of each point $(x_0, 0, 1/4)$, $x_0 \in \mathbb{C}$. Namely, for every $x_0 \in \mathbb{C}$ there exists a polydisk $\Delta(x_0)$*

$$
(2.8) \quad \Delta(x_0) = \left\{ (t, v, w) \in \mathbb{C}^3 : |t - x_0| < 1 \ , \ |v| < R^{-1} \ , \ \left|w - \frac{1}{4}\right| < \epsilon \right\}
$$

*($R \geq 1$ and $\epsilon < 1/4$) where $\Gamma$ is holomorphic and where its power series in $v$ has the form*

$$
(2.9) \quad \gamma(t, v, w) = t + \sum_{k \geq 5} \gamma_k(t, w) v^k
$$

$$
(2.10) \quad \eta(t, v, w) = v + \sum_{k \geq 4} \eta_k(t, w) v^k
$$

$$
(2.11) \quad \theta(t, v, w) = w + \sum_{k \geq 4} \theta_k(t, w) v^k
$$

*(ii) $\Gamma$ is a biholomorphism of $\Delta(x_0)$ onto its image. Its inverse is holomorphic near $(x_0, 0, 1/4)$ and $\Gamma^{-1} = Id + O(v_1^4)$.*

*(iii) Let $t \to (v(t), w(t))$ be a solution of (2.6) such that $(t, v(t), w(t)) \in \Delta(x_0)$ for $t$ in some disk $|t - t_1| < \epsilon'$.*

*Then the function $t \to \gamma(t, v(t), w(t))$ is a biholomorphism of the disk $|t - t_1| < \epsilon'$ onto its image.*

The proof of Proposition 1 is given in § 3.1.

The local equivalence of (2.5) and (2.6) stated in Proposition 1 is straightforwardly translated into a local equivalence of (2.1) and (2.2):

**Corollary 2.** *There exists a map* $(t, u, u') \mapsto \tilde{\Gamma}(t, u, u') = (x, y, y')$, *in the extended phase space, of the form*

(2.12)
$$\begin{cases} x = \gamma\left(t, -\frac{u'}{u^2}, \frac{u^3}{u'^2}\right) \\ y = u \\ y' = -u^2 \eta\left(t, -\frac{u'}{u^2}, \frac{u^3}{u'^2}\right) = u'\left[1 + \left(\frac{u'}{u^2}\right)^4 \eta_1\left(t, -\frac{u'}{u^2}, \frac{u^3}{u'^2}\right)\right] \end{cases}$$

*where* $\gamma, \eta, \eta_1$ *are holomorphic on the polydisks* $\Delta(x_0)$ *($x_0 \in \mathbb{C}$) which takes equation (2.1) into (2.2).*

*The transformation* $\tilde{\Gamma}$ *is one-to-one on each domain*

(2.13) $$\tilde{\Delta}(x_0) = \left\{(t, u, u') \in \mathbb{C}^3 : |t - x_0| < 1, \left|\frac{u'}{u^2}\right| < R^{-1}, \left|\frac{u^3}{u'^2} - \frac{1}{4}\right| < \epsilon\right\}$$

*If* $u(t)$ *is a solution of (2.2) such that the set* $\mathcal{T}_u = \{(t, v(t), w(t)) ; |t - t_1| < \epsilon'\}$ *is included in* $\Delta(x_0)$, *then the equality* $y(x) = u(t)$ *defines a solution of (2.1) for* $x \in \gamma(\mathcal{T}_u)$.

**Remark 3.** *If* $p \in \mathbb{C}$ *is a pole of* $u(t)$ *then* $v(t)$ *and* $w(t)$ *have removable singularities at* $p$, *with* $v(p) = 0$ *and* $w(p) = 1/4$.

As a consequence of the equivalence of Proposition 1, the Painlevé property of (2.1) follows naturally:

**Proposition 4.** *The only singularities of the solutions of (2.1) are second order poles.*

**Remark 5.** *The following description follows from the proof of Proposition 4. If* $y(x)$ *is a solution of (2.1), then whenever* $x$ *is sufficiently close to a singularity of* $y(x)$, $(x, y(x), y'(x))$ *falls in the equivalence region of (2.1) with (2.2). In fact, the complex plane is divided into a union* $\mathcal{A}$ *of nonintersecting balls – each containing a pole of* $y(x)$, *where* $(x, y(x), y'(x)) \in \text{Domain}(\Gamma)$ *– and the complement of* $\mathcal{A}$, *which is a connected set, where* $y(x)$ *is* analytic.

## 3. Proofs

### 3.1. Proof of Proposition 1. (i) The proof proceeds in several steps.

*Step 1:* A straightforward calculation shows that $\Gamma$ (cf. (2.7)) maps (2.5) to (2.6) iff $\gamma$ satisfies the equation

(3.14) $$v(L^2 \gamma) + 6(L\gamma) - (6 + w^2 v^4 \gamma)(L\gamma)^3 = 0$$

where $L$ is the linear operator

$$L = \frac{\partial}{\partial t} + \left(\frac{2}{w} - 6\right)\frac{\partial}{\partial v} + \frac{12}{v}\left(w - \frac{1}{4}\right)\frac{\partial}{\partial w}$$

and $\eta$ is given by

$$\eta = \frac{v}{L\gamma} \tag{3.15}$$

*Step 2:* There exists a formal series solution of equation (3.14) of the form (2.9) where $\gamma_k$ have convergent power series at $w = 1/4$

$$\gamma_k(t,w) = \sum_{n \geq 0} \gamma_{k,n}(t) \left(w - \frac{1}{4}\right)^n \tag{3.16}$$

with $\gamma_{k,n}$ polynomials in $t$.

Indeed, substitution of (2.9) in (3.14) leads to the recurrence

$$P_k(w)\gamma_{k+1} = h_k(\gamma_5, ..., \gamma_k, t, w) \qquad (k \geq 4) \tag{3.17}$$

where $P_k$ is the linear operator

$$P_k(w) = 144\left(w - \frac{1}{4}\right)^2 \frac{\partial^2}{\partial w^2} + \left(w - \frac{1}{4}\right) f_k(w)\frac{\partial}{\partial w} + g_k(w) \tag{3.18}$$

with

$$f_k(w) = 12(2k+1)\left(\frac{2}{w} - 6\right) \tag{3.19}$$

$$g_k(w) = (k+1)\left[\frac{2(2k+3)}{w^2} - \frac{24(k+2)}{w} + 36(k+2)\right] \tag{3.20}$$

Furthermore, $h_k(\gamma_5, ..., \gamma_k, t, w)$ is a polynomial in $\gamma_5, ..., \gamma_k$

$$h_k(\gamma_5, ..., \gamma_k, t, w) = h_k^0(\gamma_5, ..., \gamma_k, t, w) - h_k^1(\gamma_{k-7}, ..., \gamma_k, w) \tag{3.21}$$

where

$$h_k^0(\gamma_5, ..., \gamma_k, t, w) = w^2 t \delta_{k=4} - 12\partial_t \gamma_k \tag{3.22}$$
$$+ w^2 \gamma_{k-4} + 18 \sum_{i+j=k} L_i L_j + 6 \sum_{i+j+l=k} L_i L_j L_l$$
$$+ tw^2 \left\{ 3L_{k-4} + 3 \sum_{i+j=k-4} L_i L_j + \sum_{i+j+l=k-4} L_i L_j L_l \right\}$$
$$+ w^2 \left\{ 3 \sum_{i+p=k-4} \gamma_p L_i + 3 \sum_{i+j+p=k-4} \gamma_p L_i L_j \right.$$
$$\left. + \sum_{i+j+l+p=k-4} \gamma_p L_i L_j L_l \right\}$$

with the convention that the summation indices satisfy $i, j, l \geq 4$ and $p \geq 5$, where

$$(3.23) \qquad L\gamma := 1 + \sum_{k \geq 4} L_k v^k$$

$$(3.24) \qquad L_k = \partial_t \gamma_k \delta_{k \geq 5} + 12 \left( w - \frac{1}{4} \right) \partial_w \gamma_{k+1} + (k+1) \left( \frac{2}{w} - 6 \right) \gamma_{k+1}$$

and

$$h_k^1(\gamma_{k-7}, ..., \gamma_k, w) = -\frac{12}{w} \left( w - \frac{1}{4} \right) k \partial_t \gamma_k + 24 \left( w - \frac{1}{4} \right) \partial_{tw} \gamma_k + \partial_{tt} \gamma_{k-1}$$
(3.25)

We note that (3.17) gives $\gamma_{k+1}$ in terms of $\gamma_5, ..., \gamma_k$ as a solution of a second order linear inhomogeneous ODE. We are looking for functions $\gamma_{k+1}$ holomorphic at $(t, w) = (x_0, 1/4)$. The point $w = 1/4$ is a regular singular
point of (3.17) and we need to show there exist analytic solutions there.
Substituting the formal series (3.16) and the expansions of $f_k, g_k$ and $h_k$

$$f_k(w) = \sum_{n \geq 0} f_{k,n} \left( w - \frac{1}{4} \right)^n \ , \ g_k(w) = \sum_{n \geq 0} g_{k,n} \left( w - \frac{1}{4} \right)^n$$

$$h_k(\gamma_5(t, w), ..., \gamma_k(t, w), t, w) = \sum_{n \geq 0} h_{k,n}(t) \left( w - \frac{1}{4} \right)^n$$

in (3.17), we get, for $n \geq 0$,

$$(3.26) \quad [144n(n-1) + n f_{k,0} + g_{k,0}] \gamma_{k+1,n} = h_{k,n} - \sum_{p+j=n, p<n} (p f_{k,j} + g_{k,j}) \gamma_{k+1,p}$$

If for some $k$ it is true that

$$(3.27) \qquad 144n(n-1) + n f_{k,0} + g_{k,0} \neq 0$$

for all $n \geq 0$, then the system (3.26) can be solved for $\gamma_{k+1,n}$ recursively, for any values of the right side of each equation. This is the case if $k \neq 6$. Indeed,

$$(3.28) \qquad f_{k,0} = f_k(1/4) = 24(2k+1) \ , \ g_{k,0} = g_k(1/4) = 4(k+1)(k-6)$$

and the solvability conditions (3.27) reads

$$144n(n-1) + 24n(2k+1) + 4(k+1)(k-6) \neq 0$$

which holds if $k \neq 6$.

We must look at the case $k = 6$ separately. The equations for $\gamma_5, \gamma_6, \gamma_7$ are

$$P_4(w)\gamma_5 - tw^2 = 0 \tag{3.29}$$

$$P_5(w)\gamma_6 + 24\left(w - \frac{1}{4}\right)\frac{\partial^2}{\partial t \partial w}\gamma_5 + \left(\frac{20}{w} - 72\right)\frac{\partial}{\partial t}\gamma_5 = 0 \tag{3.30}$$

$$P_6(w)\gamma_7 + 24\left(w - \frac{1}{4}\right)\frac{\partial^2}{\partial t \partial w}\gamma_6 + \left(\frac{24}{w} - 84\right)\frac{\partial}{\partial t}\gamma_6 + \frac{\partial^2}{\partial t^2}\gamma_5 = 0$$

(3.31)

Direct substitution shows that the power series of $\gamma_{5,6}$ have the form

$$\gamma_5 = -\frac{t}{640} - \frac{7t}{352}\left(w - \frac{1}{4}\right) + O\left(\left(w - \frac{1}{4}\right)^2\right) \tag{3.32}$$

$$\gamma_6 = -\frac{1}{1920} + O\left(w - \frac{1}{4}\right) \tag{3.33}$$

Using (3.32) and (3.33) it follows that (3.31) also has power series series solutions $\gamma_7$ indexed by the arbitrary coefficient $\gamma_{7,0}$ (there are infinitely many such solutions because of the potential obstruction at $k = 6$, see also the note below).

By Frobenius' theory of regular singularities, the series (3.16) converge. Induction shows that the coefficients depend polynomially on the parameter $t$. In steps 3 to 8 we show that the series $\gamma_k$ (cf. (3.16)) converge in fact on a *common* polydisk.

**Note.** The special form of P1 is essential in overcoming the obstruction at $k = 6$. Generic perturbations of P1 lead to a $k = 6$ equation without solutions, implying that no (integer) power series for $\gamma$ exists.

*Step 3: A space of analytic functions.* Consider the class $\mathcal{H}_p(D)$ of holomorphic functions in the polydisk

$$D = \{(t, w) \in \mathbb{C}^2 : |t - x_0| < 1, |w - 1/4| < \epsilon\} \quad \text{(for } \epsilon \in (0, 1/4)\text{)}$$

which depend polynomially on $t$ and are continuous in the closure $\overline{D}$ of $D$. If

$$\psi(t, w) = \sum_{j=0}^{n} \psi_j(w)(t - x_0)^j \in \mathcal{H}_p(D)$$

then we define

$$|||\psi||| := \sum_{j=0}^{n} ||\psi_j(w)||, \quad \text{where } ||\psi_j|| = \sup_{|w - 1/4| < \epsilon} |\psi_j(w)| \tag{3.34}$$

We note the following estimate of $t$–derivatives

$$|||\partial_t \psi||| = \sum_{j=0}^{n} j ||\psi_j(w)|| \leq n |||\psi||| \tag{3.35}$$

and that multiplication is continuous:

$$|||\psi_1 \psi_2||| \leq |||\psi_1||| \, |||\psi_2|||$$

The usual sup norm on $\overline{D}$ is estimated in terms of the norm (3.34) by:

$$\sup_{(t,w)\in D} |\psi(t,w)| \leq |||\psi||| \tag{3.36}$$

Thus the closure $\overline{\mathcal{H}_p(D)}$ of the space $\mathcal{H}_p(D)$ in the norm $||| \cdot |||$ is contained in the space of the holomorphic functions on $D$, continuous on $\overline{D}$.

We now denote by Const. constants which are independent of $k, C, \epsilon, R$ and $x_0$.

*Step 4: For $k \geq 7$, (3.17) can be written in the form*

$$\gamma_{k+1} = N_k \gamma_{k+1} + b_k \tag{3.37}$$

*where $b_k$ depends on $\gamma_5, \ldots, \gamma_k, t$ and $w$ and with $N_k$ a linear operator satisfying*

$$N_k \mathcal{H}_p(D) \subset \mathcal{H}_p(D) \quad \text{and}$$

$$|||N_k \psi||| \leq \text{Const.}\epsilon a_\epsilon^2 |||\psi||| \tag{3.38}$$

*for all $\psi \in \mathcal{H}_p(D)$, where $a_\epsilon = (1/4 - \epsilon)^{-1}$.*

It is convenient to rewrite (3.18), (cf. (3.19), (3.20), (3.28)) as

$$P_k = \Lambda_k + S_k$$

where

$$\Lambda_k = 144 \left(w - \frac{1}{4}\right)^2 \frac{\partial^2}{\partial w^2} + 24(2k+1)\left(w - \frac{1}{4}\right)\frac{\partial}{\partial w} + 4(k+1)(k-6)$$

and

$$S_k = \left(w - \frac{1}{4}\right)[f_k(w) - f_k(1/4)]\frac{\partial}{\partial w} + [g_k(w) - g_k(1/4)]$$

Equation (3.17) becomes

$$\Lambda_k \gamma_{k+1} = h_k - S_k \gamma_{k+1} \tag{3.39}$$

The equation $\Lambda_k \psi = 0$ has two independent solutions $(w - 1/4)^{\alpha_{1,2}}$ with

$$\alpha_1 = 1 - k/6 \quad \text{and} \quad \alpha_2 = -1/6 - k/6$$

In order to write (3.39) in integral form, we interpret it as a linear inhomogeneous equation, with the inhomogeneous part $h_k - S_k \gamma_{k+1}$. Then the solutions $\gamma_{k+1}$ satisfy

$$\gamma_{k+1} = B_k(h_k - S_k \gamma_{k+1}) \tag{3.40}$$

where $B_k$ is the operator

$$B_k(\psi)(t,w) = \frac{1}{168} \int_0^1 \left(s^{-1-\alpha_1} - s^{-1-\alpha_2}\right) \psi(t, w_s)\, ds \tag{3.41}$$

where

(3.42) $$w_s = s\left(w - \frac{1}{4}\right) + \frac{1}{4}$$

The solutions of (3.39) satisfy the equation (3.37) where

(3.43) $$N_k \psi = -B_k(S_k \psi)$$

and $b_k$, as a function of $(t, w)$, is given by

(3.44) $$b_k = B_k(h_k)$$

Let $\psi \in \mathcal{H}_p(D)$. Integrating in (3.43) by parts the term containing $\partial_w \gamma_{k+1}$ we get

(3.45) $$(N_k \psi)(t, w) = -\frac{1}{168} \int_0^1 \phi_k(s, w) \psi(t, w_s) \, ds$$

where

$$\phi_k(s, w) = -\left(\alpha_1 s^{-1-\alpha_1} - \alpha_2 s^{-1-\alpha_2}\right)\left[f_k(w_s) - f_k(1/4)\right]$$
$$+ \left(s^{-1-\alpha_1} - s^{-1-\alpha_2}\right)\left[(w - 1/4)s f'_k(w_s) + g_k(w_s) - g_k(1/4)\right]$$

Direct estimates on (3.19), (3.20) give for $|w - 1/4| \leq \epsilon < 1/4$

(3.46) $$|f_k(w) - f_k(1/4)| \leq \text{Const.} k \epsilon a_\epsilon$$
(3.47) $$|f'_k(w)| \leq \text{Const.} k a_\epsilon^2$$
(3.48) $$|g_k(w) - g_k(1/4)| \leq \text{Const.} k^2 \epsilon a_\epsilon^2$$

Also,

(3.49) $$\int_0^1 |\alpha_1 s^{-1-\alpha_1} - \alpha_2 s^{-1-\alpha_2}| \, ds < \text{Const.} k^{-1}$$

Therefore, from (3.46) to (3.49) we get

(3.50) $$\int_0^1 |\phi_k(s, w)| \, ds \leq \text{Const.} \epsilon a_\epsilon^2$$

which proves (3.38).

*Step 5:* For $k \geq 7$ the functions $\partial_w \gamma_{k+1}$ satisfy an equation of the form

(3.51) $$\partial_w \gamma_{k+1} = N_k^1(\partial_w \gamma_{k+1}) + N_k^2(\gamma_{k+1}) + b_k^1$$

where $N_k^j$ are linear operators such that

(3.52) $$N_k^j \mathcal{H}_p(D) \subset \mathcal{H}_p(D) \ (j = 1, 2)$$
(3.53) $$|||N_k^1 \psi||| \leq \text{Const.} \epsilon a_\epsilon^2 |||\psi|||$$
(3.54) $$|||N_k^2 \psi||| \leq \text{Const.} a_\epsilon^3 |||\psi|||$$

for all $\psi \in \mathcal{H}_p(D)$.

Indeed, differentiation with respect to $w$ in (3.40) yields

$$(3.55) \quad \partial_w \gamma_{k+1}(t,w) = \frac{1}{168} \int_0^1 \left(s^{-\alpha_1} - s^{-\alpha_2}\right) \partial_w \left(h_k - S_k \gamma_{k+1}\right)(t,w_s) \, ds$$

and (3.51) follows with

$$(3.56) \quad (N_k^1 \psi)(t,w) = \frac{1}{168} \int_0^1 \phi_k^1(s,w) \psi(t,w_s) ds$$

where

$$\phi_k^1(s,w) = \left(\alpha_2 s^{-1-\alpha_2} - \alpha_1 s^{-1-\alpha_1}\right) \left(f_k(w_s) - f_k(1/4)\right) \\ - \left(s^{-\alpha_1} - s^{-\alpha_2}\right) \left(g_k(w_s) - g_k(1/4)\right)$$

and

$$(N_k^2 \psi)(t,w) = -\frac{1}{168} \int_0^1 \left(s^{-\alpha_1} - s^{-\alpha_2}\right) g_k'(w_s) \psi(t,w_s) ds$$

From (3.55) we get

$$(3.57) \quad b_k^1 = B_k^1 (h_k)$$

where $B_k^1$ is the linear operator

$$(3.58) \, B_k^1(\psi)(t,w) = \frac{1}{168} \left(w - \frac{1}{4}\right)^{-1} \int_0^1 \left(\alpha_1 s^{-1-\alpha_1} - \alpha_2 s^{-1-\alpha_2}\right) \psi(t,w_s) \, ds$$

The estimates (3.53) and (3.54) are straightforward.

*Step 6:* $\gamma_k(t,w)$ *is a polynomial in t, of degree at most* $(k-1)/4$. The proof is by induction on $k$. Note first that the degree in $t$ of any holomorphic solution of (3.17) cannot exceed the degree of the inhomogeneous term $h_k$. Hence (cf. (3.29)) the degree in $t$ of $\gamma_5$ is at most 1.

Assume by induction that $\deg_t \gamma_j \leq (j-1)/4$ if $5 \leq j \leq k$. Then, (cf. (3.23)) $\deg_t L_j \leq \deg_t \gamma_{j+1} \leq j/4$ for $j \leq k-1$ and we have (cf. (3.21)) $\deg_t h_k \leq k/4$.

*Step 7: Estimates of sums.* Given $\beta < 1$ there exist constants $C_1, C_2 > 0$ such that for all integers $n \geq 1$

$$(3.59) \quad n^{-\beta} < C_2 (1-\beta)(2-\beta)...(n-\beta)/n! < C_1 n^{-\beta}$$

and ($\beta = 3/2$)

$$(3.60) \quad n^{-3/2} < -C_2 (1-3/2)(2-3/2)...(n-3/2)/n! < C_1 n^{-3/2}$$

These elementary Gamma function inequalities imply that given $\beta < 1$, there is a constant $C_0 > 0$ such that for all $k$

$$\sum_{i+j=k, i,j \geq 4} i^{-\beta} j^{-\beta} \leq C_0 k^{1-2\beta} \tag{3.61}$$

$$\sum_{i+j+l=k, i,j,l \geq 4} i^{-\beta} j^{-\beta} l^{-\beta} \leq C_0 k^{2-3\beta} \tag{3.62}$$

$$\sum_{q+j=k, j \geq 4, q \geq 5} q^{-3/2} j^{-1/2} \leq C_0 k^{-1/2} \tag{3.63}$$

$$\sum_{q+i+j=k, i,j \geq 4, q \geq 5} q^{-3/2} i^{-1/2} j^{-1/2} \leq C_0 \tag{3.64}$$

$$\sum_{q+i+j+l=k, i,j,l \geq 4, q \geq 5} q^{-3/2} i^{-1/2} j^{-1/2} l^{-1/2} \leq C_0 k^{1/2} \tag{3.65}$$

*Step 8: There exist constants $\epsilon \in (0, 1/4)$, $C, R > 1$ such that $\gamma_k \in \mathcal{H}_p(D)$ and*

$$|||\gamma_j||| < C j^{-3/2} R^{j-5} \tag{3.66}$$

$$|||\partial_w \gamma_j||| < C j^{-1/2} R^{j-5} \tag{3.67}$$

*for all $j \geq 5$.*

The proof is by induction. Suppose that $\gamma_5, ..., \gamma_k \in \mathcal{H}_p(D)$ and that (3.66) and (3.67) hold for $j = 5, ..., k$; we prove these properties for $j = k+1$. The plan of the proof is as follows:

We estimate the functions: **(1)** $h_k$ (defined in (3.21), (3.22) and (3.25)); **(2)** $b_k$ (defined by (3.44)); and **(3)** $b_k^1$ (defined by (3.57)). **(4)** We obtain inequalities for $\gamma_{k+1}$ and $\partial_w \gamma_{k+1}$; and finally **(5)** we choose the constants $C, R$ and $\epsilon$.

**(1)** Estimating the terms of $h_k^0$ and $h_k^1$:

(a) From the induction hypothesis, (3.35), (3.23)) and the fact that $|w| < 1/4$, we see that

$$|||L_i||| \leq \deg_t \gamma_i \, |||\gamma_i||| + 12\epsilon |||\partial_w \gamma_{i+1}||| + (i+1)(2a_\epsilon + 6)|||\gamma_{i+1}|||$$
$$\leq \text{Const.} K_\epsilon C i^{-1/2} R^{i-4} \tag{3.68}$$

for $i \leq k-1$, where

$$K_\epsilon = 1 + a_\epsilon \tag{3.69}$$

(b) By (a) and (3.61)

$$||| \sum_{i+j=k} L_i L_j ||| \leq \text{Const.} K_\epsilon^2 \sum_{i+j=k} C^2 i^{-1/2} j^{-1/2} R^{k-8} \leq \text{Const.} K_\epsilon^2 C^2 R^{k-8}$$

(c) By (a) and (3.62)

$$||| \sum_{i+j+l=k} L_i L_j L_l ||| \leq \text{Const.} K_\epsilon^3 \sum_{i+j+l=k, i,j,l \geq 4} C^3 i^{-1/2} j^{-1/2} l^{-1/2} R^{k-12} \leq$$
$$\text{Const.} K_\epsilon^3 C^3 k^{1/2} R^{k-12}$$

We finally have

(d) $\quad |||w^2 \sum_{i+p=k-4} \gamma_p L_i||| \leq \text{Const.} K_\epsilon C^2 k^{-1/2} R^{k-13}$

(e) $\quad |||w^2 \sum_{i+j+p=k-4} \gamma_p L_i L_j||| \leq \text{Const.} K_\epsilon^2 C^3 R^{k-17}$

and

(f) $\quad |||w^2 \sum_{i+j+l+p=k-4} \gamma_p L_i L_j L_l||| \leq \text{Const.} K_\epsilon^3 C^4 k^{1/2} R^{k-21}$

Combining (a)...(f) we get

$$\text{(3.70)} \qquad |||h_k||| \leq \text{Const.} \Phi C(k+1)^{1/2} R^{k-4}$$

where, using that $R, K_\epsilon, C$ and $k$ are at least one and $\epsilon < 1/4$ we get

$$\Phi = \left(1 + |x_0|R^{-4}\right) \left(K_\epsilon R^{-1} + K_\epsilon^2 CR^{-4} + K_\epsilon^3 C^2 R^{-8} + K_\epsilon^3 C^3 R^{-12}\right)$$

**(2)** In view of (3.70)

$$\text{(3.71)} \quad |||b_k||| = |||B_k(h_k)||| \leq \frac{1}{168}|||h_k||| \int_0^1 \left(s^{-1-\alpha_1} - s^{-1-\alpha_2}\right) ds$$
$$\leq \text{Const.} \Phi C(k+1)^{-3/2} R^{k-4}$$

**(3)** The function $B_k^1(h_k)(t,w)$ is regular at $w = 1/4$ since

$$\int_0^1 \left(\alpha_1 s^{-1-\alpha_1} - \alpha_2 s^{-1-\alpha_2}\right) ds = 0$$

Its maximum over the disk $|w - 1/4| \leq \epsilon$ is then attained for $|w - 1/4| = \epsilon$ and

$$|||b_k^1||| = |||B_k^1(h_k)||| \leq \text{Const.} |||h_k||| \int_0^1 |\alpha_1 s^{-1-\alpha_1} - \alpha_2 s^{-1-\alpha_2}| ds$$
$$\leq \text{Const.} \Phi \epsilon^{-1} C(k+1)^{-1/2} R^{k-4}$$

**(4)** In this part of the proof we restrict the values of the parameters $\epsilon, R$ and $C$ (conditions 1 to 4); the consistency of the conditions is shown in part (5).
If $\lambda_\epsilon := \text{Const.} \epsilon a_\epsilon^2 < 1$ (*Condition 1*) then, using (3.38), it follows that the operator $N_k$ extends continuously on $\overline{\mathcal{H}_p(D)}$ and is a contraction there. Thus $I - N_k$ is invertible. The function $\gamma_{k+1} = (I - N_k)^{-1} b_k$ (cf. (3.37)) is analytic in $D$, continuous on $\bar{D}$. Since (by step 6) $\gamma_{k+1}$ is a polynomial in $t$, it follows that $\gamma_{k+1} \in \mathcal{H}_p(D)$.

In view of (3.38) and (3.71) it follows that

$$\text{(3.72)} \quad |||\gamma_{k+1}||| = |||(I - N_k)^{-1} b_k||| < (1 - \lambda_\epsilon)^{-1} |||b_k||| \leq$$
$$\text{Const.} \Phi(1 - \lambda_\epsilon)^{-1} C(k+1)^{-3/2} R^{k-4}$$

The induction step for $\gamma_{k+1}$ follows provided that $\text{Const.} \Phi(1 - \lambda_\epsilon)^{-1} \leq 1$ (*Condition 2*). The estimate of $\partial_w \gamma_{k+1}$ is similar: step 5 and (3.72) imply (cf. (3.56))

$$|||N_k^2 \gamma_{k+1}||| \leq \text{Const.} a_\epsilon^3 C(k+1)^{-3/2} R^{k-4}$$

If $\lambda_\epsilon^1 \equiv \text{Const}.\epsilon a_\epsilon^2 < 1$ (*Condition 3*) then (cf. step 5)

$$|||\partial_w \gamma_{k+1}||| \leq \text{Const}.(1 - \lambda_\epsilon^1)^{-1} \Phi^1 \; C(k+1)^{-1/2} R^{k-4}$$

where

$$\Phi^1 = a_\epsilon^3 k^{-1} + \epsilon^{-1} \Phi$$

The induction step follows if $\text{Const}.(1 - \lambda_\epsilon^1)^{-1} \Phi^1 \leq 1$ (*Condition 4*).

**(5)** Proving that conditions 1 through 4 can be satisfied. Let $\epsilon$ be small enough so that conditions 1 and 3 hold. It is convenient to impose $CR^{-4} \leq 1$ (*Condition 5*).

Then conditions 2 and 4 are implied by an inequality of the form

$$(3.73) \qquad C_1(\epsilon)\Big(k^{-1} + \big(R^{-1} + CR^{-4}\big)\big(1 + |x_0|R^{-4}\big)\Big) \leq 1$$

with $C_1(\epsilon)$ depending on $\epsilon$ only. If $k$ is larger than some $k_1(\epsilon)$, then $C_1(\epsilon)k^{-1} < 1/2$. $C$ can now be chosen so that the induction hypothesis ((3.66) and (3.67)) holds for $j = 5, ..., k_1(\epsilon)$. Finally, for large enough $R$, condition 5 and (3.73) are satisfied.

*Step 9: The series (2.9) converges for $(t, v, w) \in \Delta(x_0)$ (cf. (2.8))*. This is an immediate consequence of (3.66) and (3.36).

Part (ii) follows from the fact that $\Gamma$ is a small perturbation of the identity: $\Gamma = \text{Id} + v^4 \Gamma_1$ with $\Gamma_1$ holomorphic on $\Delta(x_0)$. Therefore (possibly taking a larger $R$) we have $\sup_{\Delta(x_0)} |v^4 \Gamma_1| < 1$, hence $\Gamma$ is invertible. Part (iii) follows similarly. ∎

3.2. **Proof of Proposition 4.** Let $y(x)$ be a solution of (2.1), analytic at some point $a \in \mathbb{C}$. The idea of the proof is that if $y, y'$ are not too large, regularity of (2.1) implies $y$ is analytic in a sufficiently wide region, whereas when $y, y' \to \infty$ equivalence with (2.2) applies, thus $1/y$ is locally analytic; uniform estimates of $y$ outside the equivalence regions preclude the accumulation of poles of $y$.

Consider the open disk $B \subset \mathbb{C}$ centered at $a$, of radius $r$, such that $y(x)$ is meromorphic in $B$. We show that $y(x)$ extends as a meromorphic function across $\partial B$, which implies $y(x)$ is meromorphic on $\mathbb{C}$.

Let

$$\Delta_1(x_0) = \left\{(x, v_1, w_1) \in \mathbb{C}^3 : |x - x_0| < \xi \, , \, |v_1| < S^{-1} \, , \, \left|w_1 - \frac{1}{4}\right| < \delta \right\}$$

where $\xi, S$ and $\delta$ are small enough so that the closure of $\Delta_1(x_0)$ is contained in $\Gamma(\Delta(x_0))$. Let $\tilde{\Delta}_1(x_0)$ be its representation in $(x, y, y')$.

Let $x^* \in \partial B$. Let $x_0 \in B$ such that $|x^* - x_0| < \xi$ and $y(x)$ is analytic at $x_0$.

We show that either $(x, y(x), y'(x)) \in \tilde{\Delta}_1(x_0)$, for all $x$ in a neighborhood of $x^*$, thus $y(x)$ is meromorphic at $x^*$ (cf. corollary 2) or else there is a path in $B$, ending at $x^*$, on which $y(x)$ is uniformly bounded. In the latter case, $y(x)$ is analytic at $x^*$ as follows from lemma 6.

**Lemma 6.** *Let $y(x)$ be a solution of (2.1). Let $l : [0,1] \to \mathbb{C}$ be a path, continuous on $[0,1]$, smooth on $[0,1)$, of finite length, and such that $y(x)$ is analytic at each point on $l[0,1)$ and uniformly bounded on $l[0,1)$. Then $y(x)$ is analytic at $l(1)$.*

The proof is given in §3.3.

Consider the ray $\mathcal{R}_\theta$ starting at $x_0$ through $x^*$. If $x, x' \in \mathcal{R}_\theta$ we write $x \prec x'$ when $|x - x_0| < |x' - x_0|$.

Integrating (2.1) we get

$$(3.74) \qquad (y')^2 = 4y^3 + 2xy - 2\int_{x_0}^{x} y(s)ds + C$$

Let $\epsilon_1 \in (0,1)$ be small.

In the following, Const. denotes positive constants, which may depend only on $x_0, y(x_0), y'(x_0), r, \xi, S$ and $\delta$ but not on $\epsilon_1$.

If $|y(x)| \le \epsilon_1^{-1}$ in $\mathcal{R}_\theta \cap B$ we define $l$ to be the segment $[x_0, x^*]$. Otherwise, let $x_1$ be the least $x$ in $\mathcal{R}_\theta \cap B$ with respect to $\prec$ such that $|y(x_1)| = \epsilon_1^{-1}$. We have

$$(3.75) \qquad (y')^2 = 4y^3 + 2xy - 2\int_{x_1}^{x} y(s)ds + C_1$$

where

$$|C_1| \le |C| + 2|x_0 - x_1| \sup_{x \in [x_0, x_1]} |y| \le |C| + 2\xi\epsilon_1^{-1} \le \text{Const.}\epsilon_1^{-1}$$

for small $\epsilon_1$. Let $y_1 = y(x_1)$ and $y_1' = y'(x_1)$. Then

$$(3.76) \qquad \left|\frac{y_1'^2}{4y_1^3} - 1\right| = \left|\frac{2x_1 y_1 + C_1}{4y_1^3}\right| \le \text{Const.}\epsilon_1^2$$

so that

$$\left|\frac{y_1^3}{y_1'^2} - \frac{1}{4}\right| < \text{Const.}\epsilon_1$$

Also, from (3.75),

$$\left|\frac{y_1'}{y_1^2}\right| < \text{Const.}\epsilon_1^{1/2}$$

Thus $(x_1, y_1, y_1') \in \tilde{\Delta}_1(x_0)$ if $\epsilon_1$ is small enough.

**Remark 7.** $y(x)$ is meromorphic in a neighborhood of the closure of the disk $\mathcal{D}_{x_1} = \{x; |x - x_1| < 4\epsilon_1^{1/2}\}$.

Consider the solution of (2.2), which corresponds to $y(x)$ through $\Gamma$: $u(t) = y(x)$ (cf. corollary 2), defined in a neighborhood of $t_1$, where $(t_1, u_1, u_1') = \tilde{\Gamma}^{-1}(x_1, y_1, y_1')$.

The change of variables $u = 1/Q^2$ in (2.2) gives, after one integration,

$$(3.77) \qquad (Q')^2 = 1 + K_1 Q^6$$

To estimate the constant $K_1$ we use the fact that the map $\tilde{\Gamma}^{-1}$ is close to the identity:

$$u' = y'\left[1 + \left(\frac{y'}{y^2}\right)^4 \theta_1\left(x, -\frac{y'}{y^2}, \frac{y^3}{y'^2}\right)\right]$$

with $\theta_1$ holomorphic for $(x, y, y') \in \tilde{\Gamma}\left(\tilde{\Delta}(x_0)\right)$ (thus bounded on $\tilde{\Delta}_1(x_0)$). Then

$$|K_1| = \frac{1}{4}\left|u_1'^2 - 4u_1^3\right| < \text{Const.}\epsilon_1^{-1}$$

With the substitution $Q(t) = Q(t_1) + (t - t_1) + U(t)$, (3.77) can be written in integral form as

$$(3.78) \qquad U(t) = \int_{t_1}^t \frac{K_1 Q(s)^6}{1 + \sqrt{1 + K_1 Q(s)^6}} ds := J\bigl(U\bigr)(t)$$

A straightforward calculation shows that the operator $J$ in (3.78) is a contraction in the sup norm over the closed disk $\mathcal{D}'_{t_1} = \{t : |t - t_1| \leq 8\epsilon_1^{1/2}\}$, in the ball $\|U\|_\infty \leq \epsilon_1^2$ (if $\epsilon_1$ is small).

Furthermore, $(t, u(t), u'(t)) \in \tilde{\Delta}(x_0)$ for $t \in \mathcal{D}'_{t_1}$ (and small $\epsilon_1$). Hence, from Proposition 1 (iii), $y(x) := u(t)$ is a solution of (2.1) if $t \in \mathcal{D}'_{t_1}$ (the same solution as in the beginning of § 3.2).

To estimate the domain in the $x$-plane where $y(x)$ is defined by the equivalence, we rely on the fact that $\Gamma$ is close to the identity:

$$|x - x_1| = \left|\gamma\left(t, -\frac{u'(t)}{u(t)^2}, \frac{u(t)^3}{u'(t)^2}\right) - \gamma\left(t_1, -\frac{u_1'}{u_1^2}, \frac{u_1^3}{u_1'^2}\right)\right| \geq$$

$$|t - t_1| - M_2 \max\left|\frac{u'(t)}{u(t)^2}\right|^5 \geq |t - t_1| - \text{Const.}\epsilon_1^{5/2}$$

Thus, taking $\epsilon_1$ small enough, remark 7 follows.

**Remark 8.** *$y(x)$ satisfies the estimate*

$$|y(x)| < \epsilon_1^{-1} \quad \text{for } |x - x_1| = 4\epsilon_1^{1/2}$$

*Proof.* Indeed, since $y(x) = u(t) = [Q(t_1) + (t - t_1) + U(t)]^{-2}$ (note that $Q(t_1)^{-2} = u_1$), we have $y(x) = [Q(t_1) + (x - x_1) + Y(x)]^{-2}$ where

$$|Y(x)| = \left|U(t) - \left[\frac{u'(t)}{u^2(t)}\right]^5 \gamma_1(t, v(t), w(t)) + \left[\frac{u_1}{u_1^2}\right]^5 \gamma_1(t_1, v(t_1), w(t_1))\right|$$

$$\leq \epsilon_1^2 + \text{Const.}\epsilon_1^{5/2}$$

since $\gamma_1$ is analytic on $\Delta(x_0)$, therefore bounded on $\Gamma^{-1}(\Delta_1(x_0))$.

Thus, for $x$ such that $|x - x_1| = 4\epsilon_1^{-1}$, and small $\epsilon_1$,

$$|y(x)| \leq |Q_1 + (x - x_1)|^{-2} \left(1 - \left|\frac{Y(x)}{Q_1 + (x - x_1)}\right|\right)^{-2} < \epsilon_1^{-1}$$

which proves remark 8.

If $x^* \in \mathcal{D}_{x_1}$ then, by remark 7, proposition 4 is proved. Otherwise, let $x_1', x_2$, with $x_1' \prec x_2$, be the two points of intersections of the circle $\partial \mathcal{D}_{x_1}$ with $R_\theta$. Thus $x_1' \prec x_1 \prec x_2 \prec x^*$ and $|y(x_2)| < \epsilon_1^{-1}$. The construction in step 1 is now repeated with $x_2$ instead of $x_0$ and the same $\epsilon_1$: if $|y(x)| \leq \epsilon_1^{-1}$ on $[x_2, x^*)$ then $l$ is defined as the segment $[x_2, x^*]$. Otherwise, let $x_3$ be the least point w.r.t. $\prec$ on $(x_2, x^*)$ such that $|y(x_3)| = \epsilon_1^{-1}$. As before, $y(x)$ is meromorphic in a neighborhood of $\mathcal{D}_{x_3}$ and satisfies $|y(x)| < \epsilon_1^{-1}$ on $\partial \mathcal{D}_{x_3}$.

Then, either $x^* \in \mathcal{D}_{x_3}$, in which case the Proposition is proved, or else, if $x_3' \prec x_4$ are the two points of intersection of $\partial D_{x_3}$ with $R_\theta$ then $x_1' \prec x_1 \prec x_2 \prec x_3 \prec x_4 \prec x^*$. We take the path $l$ going along $R_\theta$, from $x_0$ towards $x^*$, avoiding the disks $\mathcal{D}_{x_{2j+1}}$ by going on the upper semicircle of $\partial \mathcal{D}_{x_{2j+1}}$.

After a finite number of steps the construction stops since $|x_{2j+1} - x_{2j-1}| > |x_{2j} - x_{2j-1}| = 4\epsilon_1^{1/2}$. ∎

**Lemma 9.** *Let $y(x)$ be a solution of (2.1), which is analytic at $x_0$. Then the radius of analyticity is at least*

$$\min\left\{|y(x_0)|^{-1/2}, |y'(x_0)/2|^{-1/3}, |y(x_0)^2 + x_0/6|^{-1/4}\right\}$$

*Proof.* Straightforward estimates of Taylor series coefficients; for details see [1]. ∎

3.3. **Proof of Lemma 6.** Suppose $|y(x)| < M$ for $x \in l[0, 1)$. Let $x = l(s), s \in [0, 1)$. Equation (2.1) can be integrated once and it yields, for the solution $y(x)$

$$(3.79) \qquad y'(x)^2 = 4y(x)^3 + 2xy(x) - 2\int_{l[0,s]} y(x')\, dx' + C_0$$

Thus $\sup\{|y'(x)| : x \in l[0, 1)\} = 2M_1 < \infty$. By Lemma 9, for any $x \in l[0, 1)$, $y(x)$ is holomorphic in a disk centered at $x$ of radius at least

$$\rho(x) = \min\left\{M^{-1/2}, M_1^{-1/3}, \left(M^2 + \max_{x \in l[0,1]} |x|/6\right)^{-1/4}\right\}$$

Therefore $y(x)$ is analytic at $l(1)$. ∎

**Acknowledgments** The authors would like to thank Professor Martin Kruskal for very interesting discussions.

Ovidiu Costin, Mathematics Department, The University of Chicago, 5734 S. University Avenue, Chicago, IL 60637

*Current address*: Mathematical Sciences Research Institute, 1000 Centennial Drive, Berkeley, CA 94720-5070

*E-mail address*: `costin@math.uchicago.edu; costin@msri.org`

Rodica D. Costin, Mathematical Sciences Research Institute, 1000 Centennial Drive, Berkeley, CA 94720-5070

*E-mail address*: `rcostin@msri.org`